\documentclass [11pt,a4paper]{article}

\usepackage{amsmath}
\usepackage{amsthm}
\usepackage{graphics}
\usepackage{latexsym}
\usepackage{amssymb}
\usepackage{enumerate}

\newtheorem {lemma}{Lemma}

\newtheorem {fact}{Fact}

\newtheorem* {theorem*}{Theorem}
\newtheorem* {thm*}{Theorem}
\newtheorem* {lemma*}{Lemma}
\newtheorem* {corollary*}{Corollary}
\newtheorem* {prop*}{Proposition}
\newtheorem* {definition*}{Definition}

\def \N {\mathbb N}
\def \Z {\mathbb Z}

\def \R {\mathbb R}

\def \T {\mathbb T}
\def \ind{1\!\!1}
\def \eps {\varepsilon}

\def \err{\text{Err}}

\def\beqs{\begin{eqnarray*}}
\def\eeqs{\end{eqnarray*}}
\def\beq{\begin{eqnarray}}
\def\eeq{\end{eqnarray}}
\def\beas{\begin{eqnarray*}}
\def\eeas{\end{eqnarray*}}
\def\bea{\begin{eqnarray}}
\def\eea{\end{eqnarray}}

\def\prob        {\ensuremath{\mathbf{P}}}
\def\expect      {\ensuremath{\mathbf{E}}}
\def\var         {\ensuremath{\mathbf{Var}}}
\def\cov         {\ensuremath{\mathbf{Cov}}}
\def\osigma         {\overline \sigma}

\def\pt{\partial_t}
\def\px{\partial_x}

\def\uz{\underline{z}}
\def\uzeta{\underline{\zeta}}
\def\vareps{\varepsilon}

\def\zmin{z_{\mathrm{min}}}
\def\zmax{z_{\mathrm{max}}}
\def\thmin{\theta_{\mathrm{min}}}
\def\thmax{\theta_{\mathrm{max}}}

\title
{Between equilibrium fluctuations and Eulerian scaling:
Perturbation of equilibrium for a class of deposition models
}

\author {B\'alint T\'oth \qquad \qquad Benedek Valk\'o}

\begin{document}

\setlength{\baselineskip}{1.23\baselineskip}

\maketitle
\begin{abstract}
We investigate propagation of perturbations of equilibrium states
for a wide class of 1D interacting particle systems. The class of
systems considered incorporates zero range, $K$-exclusion,
mysanthropic, `bricklayers' models, and much more. We do not
assume attractivity of the interactions. We apply Yau's relative
entropy method rather than coupling arguments.

The result is \emph{partial extension} of T. Sepp\"al\"ainen's recent
paper \cite{seppalainen}. For $0<\beta<1/5$ fixed, we prove that,
rescaling microscopic space and time by $N$, respectively
$N^{1+\beta}$, the macroscopic evolution of perturbations of
microscopic order $N^{-\beta}$ of the equilibrium states is
governed by Burgers' equation. The same statement should hold for
$0<\beta<1/2$ as in Sepp\"al\"ainen's cited paper, but our method
does not seem to work for $\beta\ge1/5$.
\end{abstract}


\section{Introduction}
\label{intro}

In the recent paper \cite{seppalainen} T. Sepp\"al\"ainen proves
that in the so-called totally asymmetric stick process
(equivalent to Hammersley's process as seen from a traveling
second class particle),
small perturbations  of microscopic order
$N^{-\beta}$ of equilibrium states, macroscopically propagate
according to  Burgers' equation, if hydrodynamic limit is taken
where space and time are rescaled by $N$, respectively
$N^{1+\beta}$. This result is valid for any $0<\beta<1/2$ fixed
and goes even beyond the appearence of shocks in the solution of
Burgers' equation. Sepp\"al\"ainen's proof relies on the
combinatorial peculiarities of Hammersley's model and on coupling
arguments. It is conjectured in \cite{seppalainen} that the
result should be valid in much wider context, actually Burgers'
equation should govern propagation of disturbances of equilibria
(in this scaling regime) for essentially all interacting particle
systems with one conserved observable, which under Eulerian
scaling lead to a nonlinear 1-conservation law.
Sepp\"al\"ainen's cited result and also our present paper conceptually is
closely linked to the work of R. Esposito, R. Marra and H-T. Yau,
\cite{espositomarrayau}, where this kind of intermediate scaling
was first applied  for the simple exclusion model in $d=3$.

In the present paper we partially extend Sepp\"al\"ainen's
result. We prove a very similar result universally holding for a
wide class of interacting particle systems. Our proof is
structurally robust, it does not rely on any combinatorial
properties of the models considered. We apply Yau's relative
entropy method rather than coupling arguments. We pay, of course, a
price for this generality: (1) applying the relative entropy
method, our results stay valid only up to the emergence of shocks
in the Burgers' solution and (2) we can prove our theorem only
for $\beta\in(0,1/5)$ instead of the ideal  $\beta\in(0,1/2)$.

Technically speaking, the proof is a careful application of the
relative entropy method. However, we should emphasize that there
is some new idea in the `one-block replacement' step, where
the standard large deviation argument is replaced by a central
limit estimate --- and a stronger result is gotten. See Lemma
\ref{kurschak} and its proof. Also: since in our scaling regime we
have to consider \emph{mesoscopic blocks} of size $N^{2\beta}$
rather than large microscopic blocks, in the one block estimate
so-called non-gradient arguments (e.g. spectral gap estimates)
are involved.

The paper is organized as follows. In section \ref{prelim} we
present the models considered and some preliminary computations
(infinitesimal generators, equilibria, reversed processes,
eulerian hydrodynamic limits, formal perturbations). In section
\ref{result} the main result is precisely formulated in terms of
relative entropies. Section \ref{proof} contains the proof.
This is broken up in several subsections, according to what we
consider a logical structure.


\section{Preliminaries}
\label{prelim}

\subsection{The models}
\label{models}

\subsubsection{Notation, state spece}

Throughout this paper we denote by $\T^{^{_{N}}}$ the discrete tori
$\Z/N\Z$, $N\in\N$, and by $\T$ the continuous torus $\R/\Z$.

Let $\zmin,\zmax\in\Z\cup\{-\infty,\infty\}$ with  $\zmin<\zmax$, and
$S:=[\zmin,\zmax]\cap\Z$.
The state space of the interacting particles system considered is
\begin{equation*}
\Omega^{^{_{N}}}:=S^{\T^{^{_{N}}}}.
\end{equation*}
Configurations will be denoted
\begin{equation*}
\uz:=(z_j)_{j\in T^{^{_{N}}}}\in\Omega^{^{_{N}}},
\end{equation*}

\subsubsection{Rate functions, infinitesimal generator and examples}

Following
\cite{cocozza}, \cite{rezakhanlou2} and \cite{balazs2} we require that
the rate function $c:S\times S\to[0,\infty)$ satisfy the following
conditions:

\begin{enumerate}[(A)]

\item
For any $x,y\in S$
\begin{equation*}
c(\zmin,y)=0=c(x,\zmax),
\end{equation*}
Note, that this condition is restrictive only if either
$-\infty<\zmin$ or $\zmax<+\infty$. It guarantees that, with
probability 1, the local `spins' $z_j$ stay confined within the
bounds $[\zmin,\zmax]$. In order to avoid degeneracies we also
assume that for $x\in(\zmin,\zmax]$ and $y\in[\zmin,\zmax)$
\begin{equation}
\label{nondeg}
c(x,y)>0.
\end{equation}

\item
For any $x,y,z\in S$
\begin{equation*}
c(x,y)+c(y,z)+c(z,x)=c(y,x)+c(z,y)+c(x,z).
\end{equation*}

\item
For any $x,y,z\in S\setminus\{\zmin\}$
\begin{equation*}
c(x,y-1)c(y,z-1)c(z,x-1)=c(y,x-1)c(z,y-1)c(x,z-1).
\end{equation*}
This condition is equivalent to  requiring that there exist a function
$r:S\to(0,\infty)$, with $r(\zmin)=0$,
such that for any  $x,y\in S\setminus\{\zmin\}$
\begin{equation*}
\frac{c(x,y-1)}{c(y,x-1)}=\frac{r(x)}{r(y)}.
\end{equation*}
If $-\infty<\zmin$ or $\zmax<+\infty$, we formally extend $r$ to $\Z$ as
$r(x)=0$ for $x<\zmin$, and $r(x)=\infty$ for $x>\zmax$.
%

\end{enumerate}

\noindent{\bf Remarks:}
(1)
The monotonicity condition $c(x,y+1)\le c(x,y)\le c(x+1,y)$ would imply
\emph{attractivity} of the processes defined below. We do not require this
property of the rate functions. Our arguments do not rely on coupling ideas.
\\
(2)
In the case of unbounded $z$-variable, $\max\{|\zmin|,|\zmax|\}=\infty$,
we shall also impose some growth condition on the rate function $c(x,y)$. See
condition (D) below.

The elementary movements of our Markov process are:
$(z_{j},z_{j+1})$ $\to$ $(z_{j}-1,z_{j+1}+1)$
with rate $c(z_j,z_{j+1})$.
More formally, we define
$\Theta_{j}:\Omega^{^{_{N}}}\to\Omega^{^{_{N}}}$:
\begin{eqnarray*}
\big(\Theta_{j}\uz\big)_i
=
z_i-\delta_{i,j}+\delta_{i,j+1}.
\end{eqnarray*}
The infinitesimal generator of the process defined on the torus
$\T^{^{_{N}}}$ is
\begin{eqnarray*}
&&
L^{^{_{N}}}f(\uz)
=
\sum_{j\in\T^{^{_{N}}}}
c(z_j,z_{j+1})\big(f(\Theta_{j}\uz)-f(\uz)\big).
\end{eqnarray*}
Clearly, due to the nondegeneracy condition (\ref{nondeg}), the only conserved
quantity of the process is $\sum_j z_j$.

{\bf Remark on notation:}
Consequently, we shall denote by $\uz=\left(z_j\right)_{j\in\T^{^{_{N}}}}$ an
element of the state space $\Omega^{^{_{N}}}$ and by $\uzeta(s)$ the Markov
process on $\Omega^{^{_{N}}}$ with infinitesimal generator $L^{^{_{N}}}$.


There are three essentially different classes of examples.

\begin{enumerate}[(1)]

\item
\emph{Bounded occupation number.} The only example with  $\zmin=0$
and  $\zmax=1$ is the \emph{completely asymmetric  simple
exclusion model}. For any $K>0$ one can easily check that there
exists a finite-parameter family  of models with $\zmin=0$ and
$\zmax=K$ satisfying conditions A to C. These are  usually called
\emph{generalized  $K$-exclusion models}.

\item
\emph{Occupation number bounded from below.}
There exists an in\-fi\-nite-pa\-ra\-me\-ter
family of models with $\zmin=0$ and
$\zmax=+\infty$. In particular, with
\begin{equation*}
c(x,y)=r(x)=\ind_{\{x>0\}}r(x),
\end{equation*}
we get the \emph{zero range models}.

\item
\emph{Unbounded signed occupation number.} From the
infinite-parameter family of possible models with $\zmin=-\infty$
and $\zmax=+\infty$ we point out the following: let $r:\Z\to
(0,\infty)$ satisfy
\begin{equation*}
r(z)r(-z+1)=1.
\end{equation*}
%
Define
\begin{equation*}
c(x,y)=r(x)+r(-y)
\end{equation*}
Following \cite{balazs1}, \cite{balazs2} we call these models
\emph{bricklayers models}.
\end{enumerate}

If the occupation number is not bounded (i.e. the state space is
not compact) we need some additional conditions on the growth of
the rates. In order to avoid lengthy technical computations we
only consider two special cases: the \emph{zero range model} and the
\emph{bricklayers} model, defined in examples (2) and (3). For these models we
need the following extra  conditions:

\begin{enumerate}[(D)]
\item (Growth condition for zero range and bricklayers models)
\begin{enumerate}[(i)]
\item
$ \sup\limits_{x\in \N} \left|r(x+1)-r(x)\right|\le a_1<\infty. $
\item
There exists $x_0\in\N$ and $a_2>0$ such that $r(x)-r(y)\geq a_2$
for all $x\geq y+x_0$. (That means that for $x\in \N$ $r(x)$ is
essentially linear.)
\end{enumerate}
\end{enumerate}

These conditions will guarantee the existence of dynamics and
cf. \cite{landimsethuramanvaradhan} the uniform spectral gap estimate stated
in Lemma \ref{gap}.


\subsection{Equilibrium states and reversed process}
\label{equilibria}

\subsubsection{Stationary measures}

From the growth condition D it follows that
\begin{equation*}
Z:=
\sum_{n=1}^\infty\prod_{k=1}^n r(-k+1)
+
1
+
\sum_{n=1}^\infty\prod_{k=1}^n r(k)^{-1}<\infty.
\end{equation*}
We define the following probability measure on $S$
\begin{equation*}
\pi(x):=
\left\{
\begin{array}{lcr}
\displaystyle
Z^{-1}\prod_{k=1}^{x} r(k)^{-1} & \text{ if } & x\ge 0,
\\[12pt]
\displaystyle
Z^{-1}\prod_{k=1}^{-x} r(-k+1)     & \text{ if } & x\le 0.
\end{array}
\right.
\end{equation*}
For $\theta\in\R$ let
\begin{equation*}
F(\theta):=\log\sum_{z\in S} e^{\theta z}\pi(z)
\end{equation*}
and
\begin{equation*}
\thmin:=\inf\{\theta: F(\theta)<\infty\}
\qquad
\thmax:=\sup\{\theta: F(\theta)<\infty\}
\end{equation*}
For $\theta\in(\thmin,\thmax)$ we define the probability measures
\begin{equation*}
\pi_{\theta}(z):=\pi(z)\exp\{\theta z - F(\theta)\}
\end{equation*}
on $S$.
Expectation, variance and covariance with respect to the measure $\pi_\theta$
will be denoted by $\expect_\theta(\cdots)$, $\var_\theta(\cdots)$ and
$\cov_\theta(\cdots)$, respectively.

According to \cite{cocozza}, \cite{rezakhanlou2}, \cite{balazs2}, conditions
A to C guarantee that for any $\theta\in(\thmin,\thmax)$ the
product measure
\begin{equation*}
\pi^{^{_{N}}}_{\theta}:=\prod_{j\in\T^{^{_{N}}}}\pi_{\theta}.
\end{equation*}
is stationary for the Markov process. However, due to the conservation of
$\sum_j z_j$, on the finite tori $T^{^{_{N}}}$ these measures are not
ergodic. It is a standard matter to check that the measures conditioned on the
value of $\sum_j z_j$,
\begin{equation*}
\pi^{^{_{N}}}_k(\uz):=\pi^{^{_{N}}}_{\theta}(\uz|\sum_j z_j=k),
\qquad
k\in\Z\cap[N\zmin,N\zmax],
\end{equation*}
are ergodic. We shall refer to $\pi^{^{_{N}}}_{\theta}$,
respectively,  $\pi^{^{_{N}}}_k$ as \emph{grand canonical},
respectively, \emph{canonical} measures for our model. (The
different uses of the subscript should not cause any confusion.)

\subsubsection{The reversed process}

The elementary movements of the reversed process are
$( z_{j-1},z_{j})\to(z_{j-1}+1,z_{j}-1)$
with rate $c(z_{j},z_{j-1})$.

Define $\Theta_{j}^{*}:\Omega^{^{_{N}}}\to\Omega^{^{_{N}}}$,
\begin{eqnarray*}
\big(\Theta_{j}^{*}\uz\big)_i
=
z_i-\delta_{i,j}+\delta_{i,j-1}.
\end{eqnarray*}
The reversed generator on the torus $\T^{^{_{N}}}$:
\begin{eqnarray*}
L^{^{_{N}}*}f(\uz)
=
\sum_{j\in\T^{^{_{N}}}}
c(z_{j},z_{j-1})\big(f(\Theta_{j}^{*}\uz)-f(\uz)\big).
\end{eqnarray*}
Note, that the reveresed process is the same for any $\pi^{^{_{N}}}_\theta$,
$\theta\in(\thmin, \thmax)$, or $\pi^{^{_{N}}}_k$,
$k\in\Z\cap[N\zmin,N\zmax]$.

\subsubsection{Some expectations}

We denote
\begin{equation*}
v(\theta):=\expect_\theta(z)=
\sum_{z\in S}z\pi_{\theta}(z)=
F{'}(\theta).
\end{equation*}
Elementary computations show
\begin{equation*}
v{'}(\theta)=F{''}(\theta)=
\var_\theta(z)
>0,
\end{equation*}
thus $(\thmin,\thmax)\ni\theta\mapsto v(\theta)\in(\zmin,\zmax)$ is
invertible. With some abuse of notation
denote the inverse function by $\theta(v)$.

Further notation: we shall denote
\begin{eqnarray*}
\Phi_j
&:=&
c(z_{j+1},z_{j}),
\\[5pt]
\widehat{\Phi}(v)
&:=&
\expect_{\theta(v)}(\Phi_j)
=
\sum_{x,y\in S}\pi_{\theta(v)}(x)\pi_{\theta(v)}(y) c(x,y).
\end{eqnarray*}
Clearly, if $-\infty<\zmin<\zmax<\infty$ then $\widehat{\Phi}(v)$ is bounded.
On the other hand, for the zero range models and bricklayers' models with rate
function $r$ satisfying condition (D),  straightforward estimates show that
\begin{equation*}
\widehat{\Phi}(v)\leq C |v|
\end{equation*}
and also that $\Phi_j$ has finite exponential moment with
respect to any grand canonical measure.

{\bf Remark on notation of finite-base cylinder functions:} If
$\Psi:S^m\to\R$, then we shall denote
$\Psi_j:=\Psi(z_j,\dots,z_{j+m-1})$. The indices
$j\in\T^{^{_{N}}}$ are always meant periodically, mod $N$.
Expectation of $\Psi_j$ with respect to the grand canonical
measure $\pi^{^{_{N}}}_{\theta(v)}$ is denoted
\begin{equation*}
\widehat{\Psi}(v)
:=
\expect_{\theta(v)}(\Psi_j)
=
\sum_{z_1,\dots,z_m\in S}
\pi_{\theta(v)}(z_1) \dots \pi_{\theta(v)}(z_m)
\Psi(z_1,\dots,z_m).
\end{equation*}
%


\subsection{Hydrodynamic limits}
\label{hdl}

\subsubsection{Eulerian scaling and its formal perturbation}

For the local density $v(t,x)$ of the conserved quantity $\sum_jz_j$,
\emph{under Eulerian scaling},  by applying Yau's  relative entropy method
(see \cite{yau}, or chapter 6 of \cite{kipnislandim}, or section 8 of
\cite{fritz}),   one gets the pde:
\begin{equation}
\label{euler}
\pt v + \px \widehat\Phi(v)=0.
\end{equation}

\subsubsection{Perturbation of the Euler equation}

Throughout the rest of this paper $v_0\in(\zmin,\zmax)$ will be fixed and the
shorthand notation
\begin{equation}
\label{char}
a_0:=\widehat{\Phi}(v_0),
\qquad
b_0:=\widehat{\Phi}'(v_0),
\qquad
c_0:=\widehat{\Phi}''(v_0)
\end{equation}
will be used. Note that $b_0$ is the \emph{characteristic speed}
for the hyperbolic pde (\ref{euler}), corresponding to
$v_0$. Furthermore, it is assumed that  $c_0\not=0$.

We now consider a small perturbation of the trivial constant
solution $v(t,x)\equiv v_0$ of (\ref{euler}). We fix $\beta>0$ and
insert in (\ref{euler})
\begin{equation*}
v^{(\vareps)}(t,x):=
v_0
+
\vareps^{\beta}u(\vareps^{1+\beta}t,\vareps (x-b_0t)).
\end{equation*}
Letting  $\vareps\to0$, \emph{formally}  the inviscid Burgers' equation is
gotten for $u$:
\begin{equation}
\label{burgers}
\pt u + \frac{c_0}{2} \px (u^2)=0.
\end{equation}
%


\section{The main result}
\label{result}


\subsection{Further notation and terminology}
\label{terminology}

Let $v_0\in(\zmin,\zmax)$ be fixed and  $a_0$, $b_0$ and $c_0$ as defined  in
(\ref{char}), $c_0\not=0$ is assumed. We also denote $\theta_0:=\theta(v_0)$.

Furthermore, let $u(t,x)$, $t\in[0,T]$, $x\in\T$, be \emph{smooth} solution of
Burgers' equation (\ref{burgers}). We shall use as \emph{absolute reference
measure} the stationary measure
\begin{equation*}
\pi^{^{_{N}}}:=\prod_{j\in\T^{^{_{N}}}}\pi_{\theta_0}.
\end{equation*}
We define
\begin{equation*}
\theta^{^{_{N}}}(t,x)
:=
N^{\beta}
\big(
\theta\big(v_0+N^{-\beta} u(t,x-N^\beta b_0t)\big)
-
\theta_0
\big)
\end{equation*}
i.e.
$\theta( v_0 + N^{-\beta} u(t,x-N^\beta b_0t) ) = \theta_0 +  N^{-\beta}
\theta^{^{_{N}}}(t,x)$.

The partial derivatives of $\theta^{^{_{N}}}(t,x)$ are easily computed:
\begin{eqnarray}
\notag
&&
\!\!\!\!\!\!\!\!\!\!\!\!\!\!\!\!
\theta^{^{_{N}}}_x(t,x)
:=
\px\theta^{^{_{N}}}(t,x)
=
\theta'( v_0 + N^{-\beta} u(t,x-N^\beta b_0t) )
\, \px u(t,x-N^\beta b_0t)
\\[5pt]
\label{ptth}
&&
\!\!\!\!\!\!\!\!\!\!\!\!\!\!\!\!
\theta^{^{_{N}}}_t(t,x)
:=
\pt\theta^{^{_{N}}}(t,x)
=
-
\theta^{^{_{N}}}_x(t,x)
\times \big(c_0 u(t,x-N^\beta b_0t) + N^\beta b_0 \big)
\end{eqnarray}
In the computation of $\pt \theta^{^{_{N}}}$ we use the fact that $u$ is
smooth solution of (\ref{burgers}).

The \emph{time dependent reference measure} (not to be confused with the
absolute reference measure!) is
\begin{equation}
\label{nudef}
\nu^{^{_{N}}}_{t}
:=
\prod_{j\in\T^{^{_{N}}}}
\pi_{\theta_0+ N^{-\beta} \theta^{^{_{N}}}(t,j/N)}
=
\prod_{j\in\T^{^{_{N}}}}
\pi_{\theta( v_0 + N^{-\beta} u(t,j/N-N^\beta b_0t) )}.
\end{equation}
The \emph{true distribution} of our process on $\T^{^{_{N}}}$, at
macroscopic time $t$, i.e. at microscopic time $ N^{1+\beta} t$ is
\begin{equation}
\label{mudef}
\mu^{^{_{N}}}_t
:=
\mu^{^{_{N}}}_0
\exp
\big\{
N^{1+\beta} t L^{^{_{N}}}
\big\}.
\end{equation}
The Radon-Nikodym derivatives of these last two probability measures on
$\Omega^{^{_{N}}}$, with respect to the absolute reference measure
$\pi^{^{_{N}}}$,
are
\begin{eqnarray}
\label{radnikref}
f^{^{_{N}}}_t(\uz)
&:=&
\frac{d\nu^{^{_{N}}}_t}{d\pi^{^{_{N}}}}(\uz)
\\
&=&
\!\!\!\!
\prod_{j\in\T^{^{_{N}}}}
\!\!
\exp
\big\{
z_jN^{-\beta}\theta^{^{_{N}}}(t,j/N)
-F(\theta_0+N^{-\beta}\theta^{^{_{N}}}(t,j/N))
+F(\theta_0)
\big\}
\notag
\\
\notag
h^{^{_{N}}}_t(\uz) &:=&
\frac{d\mu^{^{_{N}}}_t}{d\pi^{^{_{N}}}}(\uz)=
\exp
\big\{
N^{1+\beta}tL^{^{_{N}}*}
\big\}
h^{^{_{N}}}_0(\uz)
\end{eqnarray}
%


\subsection{What is to be proved?}
\label{what}

We want to prove that  if $\mu^{^{_{N}}}_0$ is close to $\nu^{^{_{N}}}_0$, in
the sense of the relative entropy $H(\mu^{^{_{N}}}_0\,|\,\nu^{^{_{N}}}_0)$
being
small, then $\mu^{{{N}}}_t$ stays close to $\nu^{^{_{N}}}_t$ in the same sense,
uniformly for $t\in[0,T]$.

How close? Given two smooth profiles $u_i:\T\to\R$, $i=1,2$, let
\begin{equation*}
\nu^{^{_{N}}}_{i}
:=
\prod_{j\in\T^{^{_{N}}}}
\pi_{\theta( v_0 + N^{-\beta} u_i(j/N) )},
\qquad i=1,2.
\end{equation*}
Then, an easy computation shows that the relative entropy
$H(\nu^{^{_{N}}}_{2}\,|\,\nu^{^{_{N}}}_{1})$ is
\begin{eqnarray*}
H(\nu^{^{_{N}}}_{2}\,|\,\nu^{^{_{N}}}_{1})
\!\!\!\!
&=&
\!\!
\sum_{j\in\T^{^{_{N}}}}
H(
\pi_{\theta( v_0 + N^{-\beta} u_2(j/N) )}
\,|\,
\pi_{\theta( v_0 + N^{-\beta} u_1(j/N) )}
)
\notag
\\
\!\!\!\!
&=&
\!\!
N^{1-2\beta}
\theta_0'
\int_{\T}
\big( u_2-u_1 \big)
\big( u_2-\frac{F_0''\theta_0'}{2} (u_2+u_1) \big) dx
+ {\cal O}\big( N^{1-3\beta} \big),
\end{eqnarray*}
where $\theta_0':= \theta'(v_0)$ and $F_0'':=F''(\theta_0)$.
This suggests that one should prove
\begin{equation}
\label{main}
H^{^{_{N}}}(t)
:=
H(\mu^{^{_{N}}}_t\,|\,\nu^{^{_{N}}}_t)
=
o\big( N^{1-2\beta} \big),
\end{equation}
uniformly for $t\in[0,T]$.


\subsection{Main result}
\label{mainres}

Consider a generalized mysanthrope model with rate function
satisfying conditions A-D. Let $v_0\in(\zmin,\zmax)$ be fixed so that $c_0$
defined in (\ref{char}) is nonzero. Let $u:[0,T]\times\T\to\R$ be a
\emph{smooth} solution of the inviscid Burgers' equation
(\ref{burgers}). Further on, let
$\nu^{^{_{N}}}_t$, respectively, $\mu^{^{_{N}}}_t$ be the time dependent
reference measure, respectively, the true distribution of the mysanthrope
process,  defined in
(\ref{nudef}), respectively, (\ref{mudef}).

Our main result is the following

\begin{theorem*}
Let $\beta\in(0,1/5)$ be fixed.
Under the stated conditions, if
\begin{equation*}
H(\mu^{^{_{N}}}_0\,|\,\pi^{^{_{N}}})
=
{\cal O}\big( N^{1-2\beta} \big)
\end{equation*}
and (\ref{main}) holds for $t=0$, than  (\ref{main}) will hold uniformly for
$t\in[0,T]$.
\end{theorem*}

\noindent
{\sl Remark:}
The statement should hold for $\beta<1/2$, but, with our method, seemingly only
$\beta<1/5$ can be treated.

From this theorem, by applying the entropy inequality the next corollary
follows:

\begin{corollary*}
Under the conditions of the Theorem, for any smooth test function
$\varphi:\T\to\R$
\begin{equation*}
N^{-1+\beta}\sum_{j\in\T^{^{_{N}}}}
\varphi\big((j-N^{1+\beta}b_0t)/N\big)
\left(\zeta_j\big(N^{1+\beta}t\big)-v_0\right)
{\buildrel\prob\over\longrightarrow}
\int_\T\varphi(x)u(t,x)\,dx,
\end{equation*}
\vskip-10pt
\noindent
as $N\to\infty$.
\end{corollary*}


\section{Proof}
\label{proof}

Our strategy is to get a Gromwall type estimate. We shall prove
\begin{equation}
\label{gromwall}
H^{^{_{N}}}(t)- H^{^{_{N}}}(0)
\le
C \int_0^tH^{^{_{N}}}(s)ds + \err^{^{_{N}}}(t).
\end{equation}
It is assumed that $H^{^{_{N}}}(0)=o\big( N^{1-2\beta} \big)$
and the error estimate
$\err^{^{_{N}}}(t)=o\big( N^{1-2\beta} \big)$ is the main point.

\noindent
{\bf Important remark on further notation:}
In the remaining part of the paper, without loss of generality, we assume
\begin{equation*}
v_0=0,
\qquad
\theta_0=0,
\qquad
a_0=0.
\end{equation*}
This means that from now on
$z$, $v$, $\theta$, $\Phi$ and $\widehat\Phi$
stand for
$z-v_0$, $v-v_0$, $\theta-\theta_0$, $\Phi-a_0$ and $\widehat\Phi-a_0$


\subsection{Estimating $\pt H^{^{_{N}}}(t) $}
\label{dtH}

In order to prove an inequality like (\ref{gromwall}) we need to
estimate $\pt H^{^{_{N}}}(t) $.
Using the well known inequality
\[
f L \log f \leq L f
\]
which holds for every $f\geq 0$, straightforward computations lead to
\begin{equation}
\label{diffgromwallestimate}
\pt H^{^{_{N}}}(t)
\le
N^{1+\beta} \int_{\Omega^{^{_{N}}}}
\frac{L^{^{_{N}}*}f^{^{_{N}}}_t}{f^{^{_{N}}}_t} \,d\mu^{^{_{N}}}_t
-
\int_{\Omega^{^{_{N}}}} \frac{\pt f^{^{_{N}}}_t}{f^{^{_{N}}}_t}
\,d\mu^{^{_{N}}}_t.
\end{equation}
(See chapter 6 of \cite{kipnislandim} or the paper \cite{yau} for details.)

\noindent
{\bf Further remarks on notation:}
In  subsections \ref{dtH} and \ref{blockrepl} $t\in[0,T]$ will be fixed.
In order to avoid heavy notations, in these subsections we do not denote
explicitly dependence on $t$. In particular  we shall use the following
shorthand notations
\begin{eqnarray*}
&& \theta^{^{_{N}}}(x):=\theta^{^{_{N}}}(t,x), \quad
\theta^{^{_{N}}}_x(x):=\theta^{^{_{N}}}_x(t,x), \quad
\theta^{^{_{N}}}_t(x):=\theta^{^{_{N}}}_t(t,x),
\\
&&
u^{^{_{N}}}(x):=u(t,x-N^\beta b_0 t)
\end{eqnarray*}
Discrete gradient of functions $g:\T\to \R$ will be denoted
\begin{equation*}
\nabla^{^{_{N}}}g(x):=N\big(g(x+1/N)-g(x)\big).
\end{equation*}

\subsubsection{Computation of $L^{^{_{N}}*}f^{^{_{N}}}_t/f^{^{_{N}}}_t$}

After straightforward calculations we have
\begin{eqnarray*}
\frac{L^{^{_{N}}*}f^{^{_{N}}}_t}{f^{^{_{N}}}_t}(\uz)
&=&
\sum_{j\in\T^{^{_{N}}}}
\big(
e^{-N^{-1-\beta} (\nabla^{^{_{N}}}\theta^{^{_{N}}})(j/N)}
-1
\big)
\Phi_j
\\
&=&
-N^{-1-\beta}
\sum_{j\in\T^{^{_{N}}}}
\theta^{^{_{N}}}_x(j/N)
\big(
\Phi_j
-
\widehat\Phi(N^{-\beta}u^{^{_{N}}}(j/N))
\big)
\\
&&
-N^{-1-\beta}
\sum_{j\in\T^{^{_{N}}}}
\theta^{^{_{N}}}_x(j/N)
\widehat\Phi(N^{-\beta}u^{^{_{N}}}(j/N))
\\
&&
+
\sum_{j\in\T^{^{_{N}}}}
\big(
e^{-N^{-1-\beta} (\nabla^{^{_{N}}}\theta^{^{_{N}}})(j/N)}
-
e^{-N^{-1-\beta} \theta^{^{_{N}}}_x(j/N)}
\big)
\Phi_j
\\
&&
+
\sum_{j\in\T^{^{_{N}}}}
A\big( N^{-1-\beta} \theta^{^{_{N}}}_x(j/N)\big)
\Phi_j
\end{eqnarray*}
where in the last line the shorthand notation $A(x):=e^{-x}-1+x$ is used.
The main term is the first sum on the right hand side.
We introduce
\begin{eqnarray*}
\Psi_j
&:=&
\Phi_j - b_0 z_j
\\
\widehat \Psi(v)
&:=&
\expect_{\theta(v)} \big(\Psi_j\big)
=
\widehat\Phi(v)  - b_0 v
\end{eqnarray*}
and write in the main term
\begin{eqnarray*}
\Phi_j
-
\widehat\Phi(N^{-\beta}u)
=
\big(
\Psi_j
-
\widehat\Psi(N^{-\beta}u)
\big)
+
b_0
\big(
z_j
-
N^{-\beta}u
\big)
\end{eqnarray*}
Thus, eventually we get
\begin{eqnarray}
\label{firstterm}
&&
N^{1+\beta}
\int_{\Omega^{^{_{N}}}}
\frac{L^{^{_{N}}*}f^{^{_{N}}}_t}{f^{^{_{N}}}_t}
\,d{\mu}^{^{_{N}}}_t
=
\\[5pt]
\notag
&&
\phantom{MMM}
-
\sum_{j\in\T^{^{_{N}}}}
\theta^{^{_{N}}}_x(j/N)
\int_{\Omega^{^{_{N}}}}
\big(
\Psi_j
-
\widehat\Psi(N^{-\beta}u^{^{_{N}}}(j/N))
\big)
\,d{\mu}^{^{_{N}}}_t
\\[5pt]
\notag
&&
\phantom{MMM}
-
b_0
\sum_{j\in\T^{^{_{N}}}}
\theta^{^{_{N}}}_x(j/N)
\int_{\Omega^{^{_{N}}}}
\big(
z_j
-
N^{-\beta}u^{^{_{N}}}(j/N)
\big)
\,d{\mu}^{^{_{N}}}_t
\\[5pt]
\notag
&&
\phantom{MMM}
+
\err_1^{^{_{N}}}(t)+\err_2^{^{_{N}}}(t)+\err_3^{^{_{N}}}(t),
\notag
\end{eqnarray}
where the error terms are
\begin{eqnarray}
\label{err1}
&&
\!\!\!\!\!\!\!\!\!\!\!\!\!\!\!\!\!\!\!
\err_1^{^{_{N}}}(t)
=
-
\sum_{j\in\T^{^{_{N}}}}
\theta^{^{_{N}}}_x(j/N)
\widehat\Phi(N^{-\beta}u^{^{_{N}}}(j/N)),
\\
\label{err2}
&&
\!\!\!\!\!\!\!\!\!\!\!\!\!\!\!\!\!\!\!
\err_2^{^{_{N}}}(t)
=
N^{1+\beta}\!\!
\sum_{j\in\T^{^{_{N}}}}\!\!
\big(
e^{-N^{-1-\beta} (\nabla^{^{_{N}}}\theta^{^{_{N}}})(j/N)}
\!
-
\!
e^{-N^{-1-\beta} \theta^{^{_{N}}}_x(j/N)}
\big)
\!\!
\int_{\Omega^{^{_{N}}}}
\!\!\!\!\!\!
\Phi_j
\,d{\mu}^{^{_{N}}}_t
\!\!
,
\\
\label{err3}
&&
\!\!\!\!\!\!\!\!\!\!\!\!\!\!\!\!\!\!\!
\err_3^{^{_{N}}}(t)
=
N^{1+\beta}
\sum_{j\in\T^{^{_{N}}}}
A
\big(
N^{-1-\beta} \theta^{^{_{N}}}_x(j/N)
\big)
\int_{\Omega^{^{_{N}}}}
\Phi_j
\,d{\mu}^{^{_{N}}}_t.
\end{eqnarray}

\subsubsection{Computation of $\pt f^{^{_{N}}}_t/f^{^{_{N}}}_t$}

Now we turn our attention to the second term on the right side of
(\ref{diffgromwallestimate}). From (\ref{radnikref}) and (\ref{ptth}) we
get:
\begin{eqnarray*}
\frac{\pt f^{^{_{N}}}_t}{f^{^{_{N}}}_t}(\uz)
&=&
\!\!\!
\phantom{-}
\sum_{j\in\T^{^{_{N}}}}
N^{-\beta}\theta^{^{_{N}}}_t(j/N)
\big(
z_j-N^{-\beta}u^{^{_{N}}}(j/N)
\big)
\\
&=&
\!\!\!
-
\sum_{j\in\T^{^{_{N}}}}
\theta^{^{_{N}}}_x(j/N)
\big(
c_0 N^{-\beta}u^{^{_{N}}}(j/N) + b_0
\big)
\big(
z_j-N^{-\beta}u^{^{_{N}}}(j/N)
\big)
\end{eqnarray*}
In the last sum we write
\begin{equation*}
c_0 N^{-\beta}u
=
\widehat\Psi'(N^{-\beta}u)
-
\big(
\widehat\Psi'(N^{-\beta}u)
-
c_0 N^{-\beta}u
\big)
\end{equation*}
and note that the second term is a small error.

Eventually we get:
\begin{eqnarray}
\label{secondterm}
&&
\!\!\!\!\!\!\!\!
-\int_{\Omega^{^{_{N}}}}
\frac{\pt f^{^{_{N}}}_t}{f^{^{_{N}}}_t}
\,d{\mu}^{^{_{N}}}_t(\uz)
=
\\[5pt]
\notag
&&
\phantom{MI}
\phantom{-}
\sum_{j\in\T^{^{_{N}}}}
\theta^{^{_{N}}}_x(j/N)
\widehat
\Psi'(N^{-\beta} u^{^{_{N}}}(j/N))
\int_{\Omega^{^{_{N}}}}
\big(
z_j-N^{-\beta}u^{^{_{N}}}(j/N)
\big)
\,d{\mu}^{^{_{N}}}_t(\uz)
\\[5pt]
\notag
&&
\phantom{MI}
-
b_0
\sum_{j\in\T^{^{_{N}}}}
\theta^{^{_{N}}}_x(j/N)
\int_{\Omega^{^{_{N}}}}
\big(
z_j
-
N^{-\beta}u^{^{_{N}}}(j/N)
\big)
\,d{\mu}^{^{_{N}}}_t
+
\err_4^{^{_{N}}}(t)
\notag
\end{eqnarray}
where
\begin{eqnarray}
\label{err4}
&&
\!\!\!\!\!\!\!\!\!
\err_4^{^{_{N}}}(t)
=
-
\sum_{j\in\T^{^{_{N}}}}
\theta^{^{_{N}}}_x(j/N)
\big(
\widehat \Psi'(N^{-\beta} u^{^{_{N}}}(j/N)
-
c_0 N^{-\beta} u^{^{_{N}}}(j/N)
\big)
\\[5pt]
\notag
&&
\phantom{MMMMMMMMM}
\times
\int_{\Omega^{^{_{N}}}}
\big(
z_j-N^{-\beta}u^{^{_{N}}}(j/N)
\big)
\,d{\mu}^{^{_{N}}}_t(\uz).
\end{eqnarray}
Note that, when inserting in (\ref{diffgromwallestimate}), the
second sums on
the right hand side of (\ref{firstterm}) and (\ref{secondterm})
cancel out.

\subsubsection{Blocks}

Throughout the paper the one-block size $l$ will be chosen, depending on the
system size $N$,  so that asymptotically
\begin{equation*}
l\gg N^{2\beta}.
\end{equation*}
We introduce  the block averages
\begin{equation*}
\Psi^l_j
:=
l^{-1}
\sum_{i=0}^{l-1}\Psi_{j+i},
\qquad
z^l_j
:=
l^{-1}\sum_{i=0}^{l-1}z_{j+i}.
\end{equation*}
The main terms (i.e. the first sums on the right hand side)
in (\ref{firstterm}), respectively, in (\ref{secondterm}) become
\begin{eqnarray}
\label{firstterm2}
\!\!\!\!\!\!
-
\sum_{j\in\T^{^{_{N}}}}
\theta^{^{_{N}}}_x(j/N)
\int_{\Omega^{^{_{N}}}}
\big(
\Psi^l_j
-
\widehat\Psi(N^{-\beta}u^{^{_{N}}}(j/N))
\big)
\,d{\mu}^{^{_{N}}}_t
+
\err_5^{^{_{N,l}}}(t),
\end{eqnarray}
respectively,
\begin{eqnarray}
\label{secondterm2}
&&
\!\!\!\!\!\!\!\!\!\!\!\!\!\!\!\!\!\!\!\!
\phantom{-}
\sum_{j\in\T^{^{_{N}}}}
\theta^{^{_{N}}}_x(j/N)
\widehat
\Psi'(N^{-\beta} u^{^{_{N}}}(j/N))
\int_{\Omega^{^{_{N}}}}
\big(
z^l_j-N^{-\beta}u^{^{_{N}}}(j/N)
\big)
\,d{\mu}^{^{_{N}}}_t(\uz)
\\[5pt]
\notag
&&
\phantom{MMMMMMMMMMM}
+
\err_6^{^{_{N,l}}}(t).
\end{eqnarray}
After rearrangement of sums the error terms
$\err_5^{^{_{N,l}}}(t)$, respectively, $\err_6^{^{_{N,l}}}(t)$ are
written as
\begin{eqnarray}
\label{err5}
&&
\err_5^{^{_{N,l}}}(t)
=
-
\sum_{j\in\T^{^{_{N}}}}
\Big(
l^{-1}\sum_{i=0}^{l-1}
\theta^{^{_{N}}}_x((j-i)/N)
-
\theta^{^{_{N}}}_x(j/N)
\Big)
\\
\notag
&&
\phantom{MMMMMMMMMMMMMMMMM}
\times
\int_{\Omega^{^{_{N}}}}
\Psi_j
\,d{\mu}^{^{_{N}}}_t(\uz)
\\[5pt]
\label{err6}
&&
\err_6^{^{_{N,l}}}(t)
=
\sum_{j\in\T^{^{_{N}}}}
\Big(
l^{-1}\sum_{i=0}^{l-1}
\theta^{^{_{N}}}_x((j-i)/N)
\widehat
\Psi'(N^{-\beta} u^{^{_{N}}}((j-i)/N))
\\
\notag
&&
\phantom{\sum_{i=0}^{l-1} MMMMMMM}
-
\theta^{^{_{N}}}_x(j/N)
\widehat
\Psi'(N^{-\beta} u^{^{_{N}}}(j/N))
\Big)
\int_{\Omega^{^{_{N}}}}
z_j
\,d{\mu}^{^{_{N}}}_t(\uz).
\end{eqnarray}

\subsubsection{Sumup and estimate of the error terms (so far)}

Summing up, from (\ref{diffgromwallestimate}), (\ref{firstterm}),
(\ref{secondterm}),  (\ref{firstterm2}) and (\ref{secondterm2}), so far we
have got:
\begin{eqnarray}
\label{sumup1}
&&
\!\!\!\!\!\!\!\!\!\!\!\!\!
\pt H^{^{_{N}}}(t)
\le
-
\sum_{j\in\T^{^{_{N}}}}
\theta^{^{_{N}}}_x(j/N)
\int_{\Omega^{^{_{N}}}}
\Big(
\Psi^l_j
-
\widehat\Psi(N^{-\beta}u^{^{_{N}}}(j/N))
\\[5pt]
\notag
&&
\phantom{MMMMMMM}
+
\widehat
\Psi'(N^{-\beta} u^{^{_{N}}}(j/N))
\big( z^l_j-N^{-\beta}u^{^{_{N}}}(j/N) \big)
\Big)
\,d{\mu}^{^{_{N}}}_t(\uz)
\\[5pt]
\notag
&&
\phantom{MMM}
+
\err_1^{^{_{N}}}(t)
+
\err_2^{^{_{N}}}(t)
+
\err_3^{^{_{N}}}(t)
\\[5pt]
\notag
&&
\phantom{MMM}
+
\err_4^{^{_{N}}}(t)
+
\err_5^{^{_{N,l}}}(t)
+
\err_6^{^{_{N,l}}}(t)
\end{eqnarray}
with the error terms given in (\ref{err1}), (\ref{err2}),
(\ref{err3}), (\ref{err4}), (\ref{err5}), (\ref{err6}), respectively.

For the estimate of the these terms we use the following lemma:

\begin{lemma}
\label{entropyestimate1} Let $\Psi:\Z^m\to \R$ be a finite
cylinder function and denote $\Psi_j:=\Psi(z_j, \dots,z_{j+m-1})$.
Assume that,  for $|\gamma|<\gamma_0$, $\expect_{\pi}(\exp\{\gamma
\Psi\})<\infty$. Then there exists a constant $C<\infty$ depending
only on $m$ and $\gamma_0$ , such that for any
$\psi_N:\T^{^{_{N}}}\to\R$,
\begin{equation*}
\sum_{j\in\T^{^{_{N}}}}
\psi_N(j)
\int_{\Omega^{^{_{N}}}}
\Psi_j
\,d{\mu}^{^{_{N}}}_t
\le
C
\max_{j\in\T^{^{_{N}}}}
\left|\psi_N(j)\right|
\big(
N^{1-\beta} + N \expect_{\pi} (\Psi)
\big),
\end{equation*}
uniformly for $t\in[0,T]$.
\end{lemma}

\begin{proof}
We may assume that
$\max_{j\in\T^{^{_{N}}}}\left|\psi_N(j)\right|=1$ and
$\expect_{\pi} \Psi(\zeta)=0$. We set $\gamma_1:=\gamma_0
N^{-\beta}<\gamma_0$ then with the entropy inequality:
\begin{eqnarray*}
&&
\!\!\!\!\!\!\!\!
\Big|
\sum_{j\in\T^{^{_{N}}}}
\psi_N(j)  \int_{\Omega^{^{_{N}}}} \Psi_j \,d{\mu}^{^{_{N}}}_t
\Big|
\\[5pt]
&&
\phantom{MMMMMM}
\le
\frac1{\gamma_1}
H({\mu}^{^{_{N}}}_t|{\pi}^{^{_{N}}})
+
\frac1{\gamma_1}
\log \expect_{\pi}
\exp
\big\{
\gamma_1
\sum_{j\in\T^{^{_{N}}}} \psi_N(j) \Psi_j
\big\}
.
\end{eqnarray*}
Applying the H\"older inequality to the second term, and using that
$\Psi_j$ and $\Psi_k$ are independent if $\left|j-k\right|>m$ we
have
\begin{equation*}
\Big|
\sum_{j\in\T^{^{_{N}}}}
\psi_N(j) \int_{\Omega^{^{_{N}}}} \Psi_j \,d{\mu}^{^{_{N}}}_t
\Big|
\le
\frac1{\gamma_1}
H({\mu}^{^{_{N}}}_t|{\pi}^{^{_{N}}})
+
\frac1{\gamma_1 m}
\sum_{j\in\T^{^{_{N}}}}
\Lambda\big(\gamma_1 m \psi_N(j)\big),
\end{equation*}
where we use the notation
$\Lambda(\gamma):=\log\expect_{\pi}\exp\{\gamma \Psi(\zeta)\}$.

Now,
$\Lambda(0)=\Lambda^\prime(0)=0$,
thus we have the  asymptotics
$\Lambda(\gamma)={\cal O}\big( \gamma^2 \big)$ for $\|\gamma\|\ll 1$.
Since $\max_{j\in\T^{^{_{N}}}}\left|\psi_N(j)\right|=1$ and
$\gamma_1={\cal O} \big( N^{-\beta} \big)$ there exists a positive constant
$C_1$ such that $\Lambda(\gamma_1 m \psi_N(j))\leq C_1 \gamma_1^2$
for every ${j\in\T^{^{_{N}}}}$. There also exists a constant $C_2$
with $H({\mu}^{^{_{N}}}_t|{\pi}^{^{_{N}}})\leq C_2 N^{1-2\beta}$.
From these the lemma follows with $C=C_2/\gamma_0+C_1 \gamma_0 m$.
\end{proof}

By Lemma \ref{entropyestimate1} and the smoothness of $u(t,x)$ we
readily get:
\begin{eqnarray*}
\label{err1order}
&& \err_{1}^{^{_{N}}}(t)
=
{\cal O}\big( N^{1-3\beta} \big),
\\
\label{err2order}
&&
\err_{2}^{^{_{N}}}(t)
=
{\cal O}\big( N^{-\beta} \big),
\\
\label{err3order}
&&
\err_{3}^{^{_{N}}}(t)
=
{\cal O}\big( N^{1-4\beta} \big),
\\
\label{err4order}
&&
\err_{4}^{^{_{N}}}(t)
=
{\cal O}\big( N^{1-4\beta} \big),
\\
\label{err5order}
&&
\err_{5}^{^{_{N,l}}}(t)
=
{\cal O}\big( N^{-\beta}l \big),
\\
\label{err6order}
&&
\err_{6}^{^{_{N,l}}}(t)
=
{\cal O}\big( N^{-2\beta} l \big).
\end{eqnarray*}

\subsection{One block replacement}
\label{blockrepl}

On the right hand side of (\ref{sumup1}) we replace the block
average $\Psi^l_j(\uz)$ by its `local equilibrium value':
$\widehat\Psi(z^l_j)$. We denote
\begin{equation}
\label{Rdef}
R(x,y):=
\widehat\Psi(x)
-
\widehat\Psi(y)
-
\widehat\Psi'(y) (x-y)
\end{equation}
Then:
\begin{eqnarray}
\notag
\pt H^{^{_{N}}}(t)
&\le&
-
\sum_{j\in\T^{^{_{N}}}}
\theta^{^{_N}}_x(j/N)
\int_{\Omega^{^{_{N}}}}
R(z^l_j,N^{-\beta}u^{^{_{N}}}(j/N))
\,d{\mu}^{^{_{N}}}_t(\uz)
\\[5pt]
\notag
&&
+
M^{^{_{N,l}}}(t)
+
{\cal O}
\big(
N^{1-3\beta}
\vee
N^{-\beta} l
\big),
\\
&\le&
\notag
\sup_{ \substack {0<t<T \\ j\in \T^N}}
\left|
\theta^{^{_N}}_x (j/N)
\right|
\sum_{j\in\T^{^{_{N}}}}
\int_{\Omega^{^{_{N}}}}
|R(z^l_j,N^{-\beta}u^{^{_{N}}}(j/N))|
\,d{\mu}^{^{_{N}}}_t(\uz)
\\
\label{blockreplacement}
&&
+
M^{^{_{N,l}}}(t)
+
{\cal O}
\big(
N^{1-4\beta}
\vee N^{-\beta} l
\big),
\end{eqnarray}
where
\begin{equation}
\label{oneblockerror}
M^{^{_{N,l}}}(t)
:=
-
\sum_{j\in\T^{^{_{N}}}}
\theta^{^{_N}}_x (j/N)
\int_{\Omega^{^{_{N}}}}
\big(
\Psi^l_j- \widehat\Psi(z^l_j)
\big)
\,d{\mu}^{^{_{N}}}_t(\uz).
\end{equation}
The estimate of $\int_0^t M^{^{_{N,l}}}(s)ds$ is done in the next
subsection, by the so-called `one block estimate'.

We estimate now the first term on the right hand side of
(\ref{blockreplacement}). Assume $N=Ml$. By the entropy inequality
\begin{eqnarray}
\label{entropyestimate2}
&&
\!\!\!\!\!\!\!\!\!\!
 \sum_{j\in\T^{^{_{N}}}}
\int_{\Omega^{^{_{N}}}}
|R(z^l_j,N^{-\beta}u^{^{_{N}}}(j/N))|
\,d{\mu}^{^{_{N}}}_t
\le
\frac{1}{\gamma}
H(\mu^{^{_{N}}}_t\,|\,\nu^{^{_{N}}}_t)
\\
\notag
&&
\phantom{MMMM}
+
\frac{1}{\gamma}
\log\Big(
\int_{\Omega^{^{_{N}}}}
\exp
\big\{
\gamma
\sum_{j\in\T^{^{_{N}}}}
|R(z^l_j,N^{-\beta}u^{^{_{N}}}(j/N))|
\big\}
\,d\nu^{^{_{N}}}_t(\uz)
\Big)
\end{eqnarray}
We estimate the integral in the second term on the right hand side
of (\ref{entropyestimate2}) using again the H\"older inequality:
\begin{eqnarray}
\notag
&&
\!\!\!\!\!\!\!\!
\int_{\Omega^{^{_{N}}}}
\exp
\big\{
\gamma
\sum_{j\in\T^{^{_{N}}}}
|R(z^l_j,N^{-\beta}u^{^{_{N}}}(j/N))|
\big\}
\,d\nu^{^{_{N}}}_t(\uz)
\\
\notag
&&
\quad
=
\int_{\Omega^{^{_{N}}}}
\exp
\big\{
\gamma
\sum_{i=1}^l\sum_{k=0}^{M-1}
|R(z^l_{kl+i},N^{-\beta}u^{^{_{N}}}((kl+i)/N))|
\big\}
\,d\nu^{^{_{N}}}_t(\uz)
\\
\notag
&&
\quad
\le
\Big(
\prod_{i=1}^l
\int_{\Omega^{^{_{N}}}}
\exp
\big\{
l\gamma
\sum_{k=0}^{M-1}
|R(z^l_{kl+i},N^{-\beta}u^{^{_{N}}}((kl+i)/N))|
\big\}
\,d\nu^{^{_{N}}}_t(\uz)
\Big)^{1/l}
\\
\label{holder}
&&
\quad
=
\Big(
\prod_{j\in\T^{^{_{N}}}}
\int_{\Omega^{^{_{N}}}}
\exp
\big\{
l\gamma
|R(z^l_{j},N^{-\beta}u^{^{_{N}}}(j/N))|
\big\}
\,d\nu^{^{_{N}}}_t(\uz)
\Big)^{1/l}
\end{eqnarray}
In the last setp we use the fact that for any fixed $i\in[1,l]$ the block
averages $\zeta^l_{kl+i}$, $k=0,1,\dots,M-1$, are independent under the
measure $\nu^{^{_{N}}}_t$.
From (\ref{Rdef}) it is easy to see that the function
\begin{eqnarray}
\label{Rmap}
x \mapsto R(x+N^{-\beta}u^{^{_{N}}}(j/N),N^{-\beta}u^{^{_{N}}}(j/N))
\end{eqnarray}
is asymptotically quadratic if $|x|\ll1$.
If the variables $z_i\in S$ are bounded than (\ref{Rmap}) is automatically
bounded. If $S$ is unbounded, but condition D holds, than  (\ref{Rmap})
is  asymptotically linearly bounded  for  $|x|\gg1$. Thus we may use Lemma
\ref{kurschak} stated below,
and eventually from (\ref{entropyestimate2}), (\ref{holder}) we get
for $\gamma_0$ sufficiently small and $l\geq 1/\gamma_0$:
\begin{equation*}
\sum_{j\in\T^{^{_{N}}}}
\int_{\Omega^{^{_{N}}}}
|R(z^l_j,N^{-\beta}u^{^{_{N}}}(j/N))|
\,d{\mu}^{^{_{N}}}_t(\uz)
\le
\frac{1}{\gamma_0} H(\mu^{^{_{N}}}_t\,|\,\nu^{^{_{N}}}_t) + C N
l^{-1}.
\end{equation*}
Consequently, using this bound in (\ref{blockreplacement}) we find
\begin{equation}
\label{sumup3}
\pt H^{^{_{N}}}(t)
\le
C  H^{^{_{N}}}(t)
+
M^{^{_{N,l}}}(t)
+
{\cal O}
\big(
N^{1-3\beta}
\vee
N^{-\beta} l
\vee
N l^{-1}
\big),
\end{equation}
holding uniformly for $t\in[0,T]$.

\begin{lemma}
\label{kurschak}
Let $\zeta_1, \zeta_2, \ldots$
be i.~i.~d.~ random variables with zero mean.
Assume
\begin{equation}
\label{lambound}
\Lambda(\lambda)
:=
\log \expect \left(e^{\lambda \zeta_i}\right)
<
\infty.
\end{equation}
Let the smooth function  $G:\R\rightarrow\R_+$ be quadratically,
respectively, linearly bounded for $|x|\ll1$, respectively,
$|x|\gg1$, i.e., $G(x)\le C_1\big(|x|\wedge (x^2/2)\big)$, with
some finite constant $C_1$. Then there exist constants
$\gamma_0>0$ and $C<\infty$, such that for any $0<\gamma<\gamma
_0$ and  $l\geq 1/\gamma_0$
\begin{equation}
\label{kur}
\expect \exp \big\{ \gamma l
G\big( (\zeta_1+\dots+\zeta_l)/l \big)
\big\}
<
C.
\end{equation}
\end{lemma}

\noindent {\sl Remarks:} (1) It is worth comparing the statement and
proof of Lemma \ref{kurschak} with the corresponding places in
previous works applying the one-block replacement, see, e.g.,
Proposition 1.6. in Part 6. of \cite{kipnislandim}. There usually
a weaker statement ($o(l)$ instead of ${\cal O}(1)$ on the right
hand side of (\ref{kur})) is gotten by use of more sophisticated
tools (large deviation principle instead of central limit
estimate). Actually, \emph{we do need} the sharper ${\cal O}(1)$ bound.
\\
(2)
The statement is easily extended: imposing  more restrictive conditions on
$\Lambda(\lambda)$, the growth condition on $G(x)$ can be relaxed. E.g.,
assuming $\Lambda(\lambda)={\cal O}(\lambda^2)$ for $|\lambda|\gg1$, we may
take $G(x)$ quadratically (rather than linearly) bounded at $|x|\gg1$.
\\
(3)
Actually,
\begin{equation*}
\lim_{l\to\infty}
\expect
\exp
\big\{
\gamma l  G\big( (\zeta_1+\dots+\zeta_l)/l \big)
\big\}
=
\left(1-\gamma \Lambda{''}(0)G{''}(0)\right)^{-1/2}.
\end{equation*}
But, since we need only the bound (\ref{kur}) and not the exact value of the
limit, we leave the proof of this as a funny exercise for the reader.

\begin{proof}
First we prove the statement with the more restrictive assumption
$\Lambda(\lambda)\leq C_2 \lambda^2/2$. Assume $\gamma<\big(C_1C_2\big)^{-1}$
and let $\xi$ be a standard Gaussian random variable, independent of the
variables $\zeta_j$. We denote by $<\dots>$ expectation with respect to the
variable $\xi$. Then we have the following chain of (in)equalities:
\begin{eqnarray*}
\expect
\exp
\big\{
\gamma l G\big((\zeta_1+\dots+\zeta_l)/l\big)
\big\}
&\leq&
\expect
\exp
\big\{
C_1\gamma \,(\zeta_1+\dots+\zeta_l )^2/(2l)
\big\}
\\
&=&
\expect
\big<
\exp
\big\{
\sqrt{C_1\gamma/l}\,(\zeta_1+\dots+\zeta_l )\,\xi
\big\}
\big>
\\
&=&
\big<
\expect
\exp
\big\{
\sqrt{C_1\gamma/l}\,(\zeta_1+\dots+\zeta_l )\,\xi
\big\}
\big>
\\
&=&
\big<
\exp
\big\{
l \Lambda \big( \sqrt{C_1\gamma/l}\,\xi\big)
\big\}
\big>
\\
&\le&
\big<
\exp
\big\{
C_2C_1\gamma\,\xi^2/2\big)
\big\}
\big>
\\
&=&
\big(1-\gamma C_1C_2\big)^{-1/2}.
\end{eqnarray*}
Now we consider the general case.
Choose $\alpha$ so large, that for any $x\in\R$
\begin{equation*}
G(x)<\ln\cosh(\alpha x).
\end{equation*}
One can do this due to the bounds imposed on $G$.
Let  $\xi_1,\xi_2,\dots $ be i.i.d  random variables
which are also independent of  the $\zeta_j$-s and have the common
distribution  $\prob\big(\xi_j=\pm\alpha\big)=1/2$.
We shall denote by $<\dots>$
expectation with respect to the random variables $\xi_j$.
We choose  $\lambda_0$, $C_3$  so that for $|\lambda|<\lambda_0$ the
quadratic bound $\Lambda(\lambda)<C_3\lambda^2/2$ holds and
fix $\gamma<\lambda_0/\alpha$.
Then we have:
\begin{eqnarray*}
&&
\expect
\exp
\big\{
\gamma l  G\big((\zeta_1+\dots+\zeta_l)/l\big)
\big\}
\leq
\cosh
\big(
\alpha\,(\zeta_1+\dots+\zeta_l)/l
\big)^{\lceil\gamma l\rceil}
\le
\\
&&
\phantom{MMMMMMMMMM}
\le
\expect
\big<
\exp
\big\{
(\xi_1+\dots+\xi_{\lceil\gamma l\rceil})
(\zeta_1+\dots+\zeta_l)/l
\big\}
\big>
\\
&&
\phantom{MMMMMMMMMM}
=
\big<
\expect
\exp
\big\{
(\xi_1+\dots+\xi_{\lceil\gamma l\rceil})
(\zeta_1+\dots+\zeta_l)/l
\big\}
\big>
\\
&&
\phantom{MMMMMMMMMM}
=
\big<
\exp
\big\{
l \Lambda
\big(
(\xi_1+\dots+\xi_{\lceil\gamma l\rceil})/l
\big)
\big\}
\big>
\\
&&
\phantom{MMMMMMMMMM}
\le
\big<
\exp
\big\{
C_3
(\xi_1+\dots+\xi_{\lceil\gamma l\rceil})^2/(2l)
\big\}
\big>.
\end{eqnarray*}
Now, since $\ln\cosh(\alpha x)\le \alpha^2 x^2/2$,
we can apply to the random  variables $\xi_j$ the argument
of the first part of this proof, with $C_2=\alpha^2$ and $C_1=C_3$, to get
\begin{eqnarray*}
\expect
\exp
\big\{
\gamma l G\big((\zeta_1+\dots+\zeta_l)/l\big)
\big\}
\leq
\big(1-\gamma C_3 \alpha^2\big)^{-1/2}.
\end{eqnarray*}
\end{proof}


\subsection{The one block estimate}
\label{1blockest}

The objective of this section is to provide an estimate for
$\int_0^tM^{^{_{N,l}}}(s)\,ds$, where $M^{^{_{N,l}}}(s)$ is given
in (\ref{oneblockerror}).

\subsubsection{Cutoff}

We cut off large values of the block averages. In case of compact state space,
i.e. $-\infty<\zmin<\zmax<\infty$ this step is completely omitted. Clearly,
\begin{equation}
\label{m=a+b}
M^{^{_{N,l}}}(t)
\leq
A^{^{_{N,l}}}_K (t)
+
B^{^{_{N,l}}}_K (t),
\end{equation}
where the terms on the right side are defined as
\begin{eqnarray*}
&&
\!\!\!\!\!\!\!\!
A^{^{_{N,l}}}_K (t)
:=
\\
\notag && \phantom{M} \sum_{j \in \T^{^{_{N}}}} \theta^{^{_N}}_x
(t,j/N) \int_{\Omega^{^{_{N}}}} \big( \Psi^l_j-\widehat\Psi(z^l_j)
\big) \ind_{\{ |z^l_j| \vee \alpha |\Psi^l_j| \le K \}}
\,d\mu^{^{_{N}}}_t(\uz),
\\ [3pt]
&&
\!\!\!\!\!\!\!\!
B^{^{_{N,l}}}_K (t)
:=
\\
\notag && \phantom{M} \sup_{\substack{0<t<T\\ j\in \T^N}} \left|
\theta^{^{_N}}_x (t,j/N)\right| \sum_{j \in \T^{^{_{N}}}}
\int_{\Omega^{^{_{N}}}} \big| \Psi^l_j-\widehat\Psi(z^l_j) \big|
\ind_{ \{|z^l_j| \vee \alpha |\Psi^l_j| > K \}} \,d
\mu^{^{_{N}}}_t(\uz),
\end{eqnarray*}
where $\alpha>0$ is a fixed constant which will only depend on
the rate function. For the estimate of $B^{^{_{N,l}}}_K (t)$ we
need the following lemma (applied with $m=1 or 2$ only):

\begin{lemma}
\label{largeKbound}
Let $\Delta:\Z^m\to \R$ be a  finite
cylinder variable. Then there exists a map $K\mapsto \epsilon(K)$, such
that $\lim_{K\to\infty}\epsilon (K)=0$ and
\begin{equation*}
\sum_{j \in \T^{^{_{N}}}}
\int_{\Omega^{^{_{N}}}}
\big|
\Psi^l_j-\widehat\Psi(z^l_j)
\big|
\ind_{\{ |\Delta^l_j|> K \}}
\,d \mu^{^{_{N}}}_t(\uz)
\le
\epsilon (K)
N^{1-2\beta}.
\end{equation*}
\end{lemma}

\begin{proof}
The entropy inequality yields:
\begin{eqnarray*}
&&
\!\!\!\!\!\!\!\!
\sum_{j \in \T^{^{_{N}}}}
\int_{\Omega^{^{_{N}}}}
\big|
\Psi^l_j-\widehat\Psi(z^l_j)
\big|
\ind_{\{|\Delta^l_j|> K\}}
\,d \mu^{^{_{N}}}_t(\uz)
\\
&&
\phantom{MM}
\le
\frac1{\gamma}
\Big(
H\big( \mu^{^{_{N}}}_t | \pi^{^{_{N}}} \big)
+
\log \expect_{\pi^{^{_{N}}}}
\exp
\big\{
\gamma \sum_{j \in \T^{^{_{N}}}}
\big|
\Psi^l_j-\widehat\Psi(\zeta^l_j)
\big|
\ind_{\{|\Delta^l_j|> K\}}
\big\}
\Big)
\end{eqnarray*}
We note  that the $j^{th}$ and $k^{th}$ terms are independent in
the last sum if $|j-k|> l+m-1$.
By the H\"older inequality, for $l\ge m$, we have
\begin{eqnarray*}
&&
\!\!\!\!\!\!\!\!
\log \expect_{\pi^{^{_{N}}}}
\exp
\big\{
\gamma \sum_{j \in \T^{^{_{N}}}}
\big|
\Psi^l_j-\widehat\Psi(\zeta^l_j)
\big|
\ind_{\{|\Delta^l_j|> K\}}
\big\}
\\
&&
\phantom{MMMMMMM}
\le
Nl^{-1}
\log \expect_{\pi^{^{_{N}}}}
\exp
\big\{
2l \gamma
\big|
\Psi^l-\widehat\Psi(\zeta^l)
\big|
\ind_{\{|\Delta^l|> K\}}
\big\}.
\end{eqnarray*}
Next we use Cauchy-Schwarz inequality:
\begin{eqnarray*}
&&
\!\!\!\!\!\!\!\!
\expect_{\pi^{^{_{N}}}}
\exp
\big\{
2l \gamma
\big|
\Psi^l-\widehat\Psi(\zeta^l)
\big|
\ind_{\{|\Delta^l|> K\}}
\big\}
\\
&&
\phantom{MMM}
\le
1
+
\expect_{\pi^{^{_{N}}}}
\left(
\ind_{\{|\Delta^l|> K\}}
\exp
\big\{
2l \gamma
\big|
\Psi^l-\widehat\Psi(\zeta^l)
\big|
\big\}
\right)
\\
&&
\phantom{MMM}
\le
1
+
\Big\{
\prob_{\pi^{^{_{N}}}}
\big(|\Delta^l|> K\big)
\Big\}^{1/2}
\Big\{
\expect_{\pi^{^{_{N}}}}
\exp
\big\{
2l \gamma
\big|
\Psi^l-\widehat\Psi(\zeta^l)
\big|
\big\}
\Big\}^{1/2}.
\end{eqnarray*}
From standard large deviation arguments it follows that there exists a
function  $[0,\infty)\ni\gamma\mapsto\Lambda(\gamma)\in[0,\infty)$ (finite for
any  finite $\gamma$!), such that
\begin{equation*}
\expect_{\pi^{^{_{N}}}}
\exp
\big\{
2l \gamma
\big|
\Psi^l-\widehat\Psi(\zeta^l)
\big|
\big\}
\le
\exp
\big\{
l\Lambda(\gamma)
\big\}.
\end{equation*}
On the other hand, using again a H\"older bound
and a standard large deviation estimate, for large $l$ we have
\begin{equation*}
\prob_{\pi^{^{_{N}}}}
\big(|\Delta^l|> K\big)
\le
m \exp \big\{ - l I(K)/(2m) \big\},
\end{equation*}
where $x\mapsto I(x)$ is the rate function
\begin{equation*}
I(x):=
\sup_{\lambda}
\left(
\lambda x
-
\log\expect_{\pi^{^{_{N}}}}
\exp
\big\{
\lambda \Delta
\big\}
\right).
\end{equation*}
We define
\begin{equation*}
\gamma(K):= \sup\{\gamma: \Lambda(\gamma)<
I(K)/(2m)\}\wedge K.
\end{equation*}
Since $\lim_{x\to\infty}I(x)=\infty$, we also have
$\lim_{K\to\infty}\gamma(K)=\infty$.
Now, putting together all our estimates, we get
\begin{eqnarray*}
&&
\!\!\!\!\!\!\!\!\!\!\!\!\!\!\!\!
\sum_{j \in \T^{^{_{N}}}}
\int_{\Omega^{^{_{N}}}}
\big|
\Psi^l_j-\widehat\Psi(z^l_j)
\big|
\ind_{\{|\Delta^l_j|> K\}}
\,d \mu^{^{_{N}}}_t(\uz)
\\
&&
\phantom{MMMMMMMMMM}
\le
\frac{1}{\gamma(K)}
\Big(
H\big( \mu^{^{_{N}}}_t | \pi^{^{_{N}}} \big)
+
Nl^{-1}(1+\sqrt m)
\Big).
\end{eqnarray*}
Noting that $H\big(\mu^{^{_{N}}}_t|\pi^{^{_{N}}}\big)={\cal
O}\big( N^{1-2\beta} \big)$ and $l\ge C N^{2\beta}$, the lemma
follows with $\epsilon(K)=C \gamma(K)^{-1}$.
\end{proof}

It is easy to see, that the functions $\Delta_j=z_j$ and
$\Delta_j=\Psi_j$ satisfy the conditions of the Lemma
\ref{largeKbound}, thus it follows that there exists a map
$K\rightarrow \epsilon(K)$ with $\lim_{K\rightarrow \infty}
\epsilon(K)=0$ and
\begin{equation}
\label{bbound}
B^{^{_{N,l}}}_K (t)
\leq
\epsilon(K) N^{1-2 \beta}.
\end{equation}

\subsubsection{General tools}

We collect in this paragraph the general, \emph{model independent facts}
used in the one-block estimate.

Let $\uzeta(s)$ be a Markov process on the countable state space $\Omega$,
with ergodic  stationary measure $\pi$.
Denote by $L$ and $L^*$  the infinitesimal generator and its adjoint, acting
on $L^2(\Omega,\pi)$.
We denote by $D(f)$ the Dirichlet form associated with the generator $L$ and
stationary measure $\pi$:
\begin{equation*}
D(f)
:=
-\int_{\Omega}fLf\,d\pi
=
-\int_{\Omega}fL^*f\,d\pi
\end{equation*}
The \emph{spectral gap} of the infinitesimal generator $L$ is $\rho^{-1}$
defined by
\begin{equation*}
\rho=\rho(L):=
\sup_{f\in L^2(\Omega,\pi)}
\frac{\var_\pi(f)}{D(f)}\in(0,\infty].
\end{equation*}
Actually, this means that  $(L+L^*)/2$, the symmetric part of $L$,
has a gap of size $\rho^{-1}$ in its spectrum, immediately to the
left of the eigenvalue $0$.

If  $V:\Omega\rightarrow \R$ is  a bounded measurable function we
denote
\begin{equation*}
\osigma \left( L+V(\cdot)\right)
:=
\sup\big\{
\mathrm{spectrum \ \ of \ \ } (L+L^*)/2 + V(\cdot)
 \big\}.
\end{equation*}

The following statement is the variational characterization of the `top of the
specrtrum' of a self-adjoint operator over a Hilbert space. It can be found in
any introductory textbook on functional analysis.

\begin{fact}
\label{varform}
For $\osigma \left( L+V(\cdot)\right)$
the following variational formula holds:
\begin{equation}
\label{varformeq}
\osigma\left(L + V(\cdot) \right)
=
\sup_{h}
\Big(
\int_{\Omega} V(\cdot)  h \,d \pi
-
D\big(\sqrt{h}\,\big)
\Big),
\end{equation}
where the supremum is taken over all probability densities with respect to the
stationary measure $\pi$.
\end{fact}

The second fact is a perturbative estimate of
$\osigma \left(L+\eps V(\cdot)\right)$. It can be found, e.g., as Theorem
1.1 in Appendix 3 of \cite{kipnislandim}.

\begin{fact}
\label{gappert}
If  $V:\Omega\to\R$ has  zero mean, i.e., $\int_{\Omega}V\,d\pi=0$,
then, for every
$\eps<\left( 2\left\| V \right\|_{\infty} \rho(L)\right)^{-1}$
\begin{equation}
\label{gappertbound}
\osigma \left(L+\eps V (\cdot)\right)
\leq
\frac{\eps^2 \rho(L)} {1-2\left\| V \right\|_{\infty} \eps \rho(L)}
\var_{\pi}(V).
\end{equation}
\end{fact}

The third general fact to be used is a direct consequence of the
Feynman-Kac formula and straighforward euclidean (inner product)
manipulations. Its proof can be found, e.g., in \cite{komoriya} or
as Lemma 7.2 in Appendix 1 of \cite{kipnislandim} .

\begin{fact}
\label{FeymanKac}
Assume now that
$V:\R_+\times\Omega\to \R$ is a bounded function.
The following bound holds
\begin{equation}
\label{fkbound}
\expect_{\pi}
\exp
\Big\{
\int_0^{t} V \left(s,\uzeta (s) \right) \,ds
\Big\}
\leq
\exp
\Big\{
\int_0^{t} \osigma \left( L+V(s,\cdot) \right)\,ds
\Big\},
\end{equation}
where now $\expect_{\pi}$ denotes expectation over the Markov chain
trajectories started from the stationary initial measure $\pi$.
\end{fact}

\subsubsection{Notations}

We shall use the notation $\mu^{^{_{N}}}$, respectively, $\mu^{^{_{l}}}$ for a
\emph{generic}  probability measure on $\Omega^{^{_{N}}}$, respectively,
$\Omega^{^{_{l}}}$. We shall denote by
$h^{^{_{N}}}(\uz)$, respectively, $h^{^{_{l}}}(\uz)$
their  Radon-Nikodym derivatives with respect to the absolute reference
measures  $\pi^{^{_{N}}}$, respectively, $\pi^{^{_{l}}}$.
Further on $\mu^{^{_{N,l,j}}}$ will denote the $[j,\dots,j+l-1]$ marginal of
$\mu^{^{_{N}}}$ and
$\mu^{^{_{N,l}}}:=N^{-1}\sum_{j\in\T^{^{_{N}}}}\mu^{^{_{N,l,j}}}$ the average
$l$-dimensional marginal of $\mu^{^{_{N}}}$. Correspondingly,
$h^{^{_{N,l,j}}}(\uz)$, respectively, $h^{^{_{N,l}}}(\uz)$ will denote the
Radon-Nikodym
derivatives of $\mu^{^{_{N,l,j}}}$, respectively, $\mu^{^{_{N,l}}}$, with
respect to the absolute reference measure $\pi^{^{_{l}}}$

For $k\in\Z$ fixed we denote:
\begin{eqnarray*}
&&
\Omega^{^{_{l}}}_k
:=
\big\{
\uz\in\Omega^{^{_{l}}}: \sum_{i=1}^l z_i=k
\big\},
\\
&&
m^{^{_{l}}}_k
:=
\pi^{^{_{l}}}(\Omega^{^{_{l}}}_k),
\quad
\\
&&
w^{^{_{l}}}_k
:=
\mu^{^{_{l}}}(\Omega^{^{_{l}}}_k),
\\
&&
\pi^{^{_{l}}}_k(\uz)
:=
\pi^{^{_{l}}}\big( \uz\, \big| \,\sum_{i=1}^l z_i=k \big)
=
\ind_{\{\uz\in\Omega^{^{_{l}}}_k\}}
\frac{\pi^{^{_{l}}}(\uz)}{m^{^{_{l}}}_k},
\\
&&
\mu^{^{_{l}}}_k(\uz)
:=
\mu^{^{_{l}}}\big( \uz\, \big| \,\sum_{i=1}^l z_i=k \big)
=
\ind_{\{\uz\in\Omega^{^{_{l}}}_k\}}
\frac{\mu^{^{_{l}}}(\uz)}{w^{^{_{l}}}_k},
\\
&&
h^{^{_{l}}}_k(\uz)
:=
\ind_{\{\uz\in\Omega^{^{_{l}}}_k\}}
\frac{\mu^{^{_{l}}}_k(\uz)}
     {\pi^{^{_{l}}}_k(\uz)}
=
\ind_{\{\uz\in\Omega^{^{_{l}}}_k\}}
\frac{m^{^{_{l}}}_k}{w^{^{_{l}}}_k} h^{^{_{l}}}(\uz).
\end{eqnarray*}

Denote by $D^{^{_{N}}}$, $D^{^{_{l}}}$ respectively
$D^{^{_{l}}}_k$ the following Dirichlet forms
\begin{eqnarray*}
D^{^{_{N}}}(f) &:=& \frac12 \sum_{i=1}^{N} \int_{\Omega^{^{_{N}}}}
c(z_i,z_{i+1}) \left( f(\Theta_i\uz)-f(\uz) \right)^2 \,d
\pi^{^{_{N}}}(\uz)
\\
&=&  \frac12 \sum_{i=1}^{N} \int_{\Omega^{^{_{N}}}} c(z_i,z_{i-1})
\left( f(\Theta_i^{*}\uz)-f(\uz) \right)^2 \,d \pi^{^{_{N}}}(\uz)
\\
D^{^{_{l}}}(f)
&:=&
 \frac12 \sum_{i=1}^{l-1}
\int_{\Omega^{^{_{l}}}}
c(z_i,z_{i+1})
\left(
f(\Theta_i\uz)-f(\uz)
\right)^2
\,d \pi^{^{_{l}}}(\uz)
\\
&=&
 \frac12 \sum_{i=2}^{l}
\int_{\Omega^{^{_{l}}}}
c(z_i,z_{i-1})
\left(
f(\Theta_i^{*}\uz)-f(\uz)
\right)^2
\,d \pi^{^{_{l}}}(\uz)
\\
D^{^{_{l}}}_k(f) &:=&  \frac12 \sum_{i=1}^{l-1}
\int_{\Omega^{^{_{l}}}_k} c(z_i,z_{i+1}) \left(
f(\Theta_i\uz)-f(\uz) \right)^2 \,d\pi^{^{_{l}}}_k(\uz)
\\
&=&
 \frac12 \sum_{i=2}^{l}
\int_{\Omega^{^{_{l}}}_k}
c(z_i,z_{i-1})
\left(
f(\Theta^{*}_i\uz)-f(\uz)
\right)^2
\,d \pi^{^{_{l}}}_k(\uz).
\end{eqnarray*}
In the definition of $D^{^{_{N}}}$ \emph{periodic}, in that of
$D^{^{_{l}}}$ and $D^{^{_{l}}}_k$ \emph{free}  boundary conditions
are understood.

It is easy to check that for any probability measure
$\mu^{^{_{l}}}$ on $\Omega^{^{_{l}}}$
\begin{equation}
\label{dirsum}
D^{^{_{l}}}
\big(
\sqrt{h^{^{_{l}}}}
\,\big)
=
\sum_{k\in\Z} w^{^{_{l}}}_k
D^{^{_{l}}}_k
\big(
\sqrt{h^{^{_{l}}}_k}
\,\big).
\end{equation}
Further on,
using convexity of the Dirichlet form one can readily prove that
\begin{equation}
\label{convex} D^{^{_{N}}} \big( \sqrt{h^{^{_{N}}}} \,\big) \geq
\frac{1}{l} \sum_{j\in \T^N} D^l \big( \sqrt{h^{^{_{N,l,j}}}}
\,\big).
\end{equation}

\subsubsection{Applying F-K formula}

We return now to the concrete computations.
Before the estimate of
\\
$\int_0^tA^{^{_{N,l}}}_K (s)ds$ we need some more notation (we do
not denote explicitly dependence on the cutoff):
\begin{eqnarray*}
&& V^{^{_{N,l}}}_j(\uz) := \big( \Psi^l_j-\widehat\Psi(z^l_j)
\big) \ind_{\{|z^l_j| \vee \alpha|\Psi^l_j| \le K\}},
\\[5pt]
&&
V^{l}(\uz)
:=
V^{^{_{N,l}}}_1(\uz),
\\[5pt]
&& V^{^{_{N,l}}}_j (t,\uz)
:= \theta^{^{_N}}_x (N^{-(1+\beta)}t,j/N)
V^{^{_{N,l}}}_j(\uz),
\\[5pt]
&&
V^{^{_{N,l}}}(t,\uz)
:=
\sum_{j\in \T^{^{_{N}}}} V^{^{_{N,l}}}_j (t,\uz).
\end{eqnarray*}
We denote by $\uzeta^{^{_{N}}}(t)$ the Markov process on
$\Omega^{^{_{N}}}$ with infinitisimal generator $L^{^{_{N}}}$ and
by $\expect_{\mu_0^{^{_{N}}}}$, respectively,
$\expect_{\pi^{^{_{N}}}}$ the \emph{path measure} of this process
starting with initial distribution $\mu_0^{^{_{N}}}$,
respectively, $\pi^{^{_{N}}}$.

By the definitions and the entropy inequality we have
\begin{eqnarray*}
&&
\!\!\!\!\!\!\!\!
\int_0^t A^{^{_{N,l}}}_K (s) \,ds
\ \ \ = \ \ \
\frac{1}{N^{1+\beta}}
\expect_{\mu^{N}_0}
\Big(
\int_0^{N^{1+\beta}t} V^{^{_{N,l}}}\big(s,\uzeta^N(s)\big)\,ds
\Big)
\\
\notag
&&
\!\!\!\!\!\!\!\!
\leq
\frac{1}{\gamma N^{1+\beta}}
\Big(
H(\mu^{^{_{N}}}_0| \pi^{^{_{N}}})
+
\log \expect_{\pi^{^{_{N}}}}
\exp
\Big\{
\int_0^{N^{1+\beta}t}
\gamma V^{^{_{N,l}}}\big(s,\uzeta^N (s)\big)\,ds
\Big\}
\Big).
\end{eqnarray*}
We apply the Feynman-Kac bound (\ref{fkbound}) and the variational formula
(\ref{varformeq}) to the second term on the right hand side of the last
inequality:
\begin{eqnarray}
\label{elso}
&&
\log \expect_{\pi^{^{_{N}}}}
\exp
\Big\{
\int_0^{N^{1+\beta}t}
\gamma V^{^{_{N,l}}}\big(s,\uzeta^{^{_{N}}} (s)\big)\,ds
\Big\}
\\
\notag
&&
\phantom{MMMM}
\le
\int_0^{N^{1+\beta}t}
\osigma
\big(
L^{^{_{N}}}+\gamma V^{^{_{N,l}}}(s,\cdot)
\big)
\,ds
\\
\notag
&&
\phantom{MMMM}
=
\int_0^{N^{1+\beta}t}
\sup_{h^{^{_{N}}}}
\Big(
\int_{\Omega^{^{_{N}}}}
\gamma V^{^{_{N,l}}}(s,\cdot)  h^{^{_{N}}}
\,d \pi^{^{_{N}}}
-
D^{^{_{N}}}
\big(\sqrt{h^{^{_{N}}}}
\,\big)
\Big)
\,ds.
\end{eqnarray}
Using (\ref{convex}) we bound the integrand in the last expression
\begin{eqnarray}
\label{masodik}
&&
\sup_{h^{^{_{N}}}}
\Big(
\int_{\Omega^{^{_{N}}}}
\gamma V^{^{_{N,l}}}(s,\cdot)  h^{^{_{N,l}}}
\,d \pi^{^{_{l}}}
-
D^{^{_{N}}}
\big(
\sqrt{h^{^{_{N}}}}
\,\big)
\Big)
\\
\notag
&&
\phantom{MMMMMM}
=
\sup_{h^{^{_{N}}}} \Big( \sum_{j\in\T^{^{_{N}}}}
\int_{\Omega^{^{_{l}}}} \gamma V^{^{_{N,l}}}_j(s,\cdot)
h^{^{_{N,l,j}}} \,d \pi^{^{_{l}}}
-
D^{^{_{N}}}
\big(
\sqrt{h^{^{_{N}}}}
\,\big)
\Big)
\\
\notag
&&
\phantom{MMMMMM}
\le
\sup_{h^{^{_{N}}}} \sum_{j\in\T^{^{_{N}}}} \Big(
\int_{\Omega^{^{_{l}}}} \gamma V^{^{_{N,l}}}_j(s,\cdot)
h^{^{_{N,l,j}}} \,d \pi^{^{_{l}}}
-
\frac{1}{l}D^{^{_{l}}} \big( \sqrt{h^{^{_{N,l,j}}}} \,\big) \Big)
\\
\notag
&&
\phantom{MMMMMM}
\le
\frac1l
\sum_{j\in\T^{^{_{N}}}}
\sup_{h^{^{_{l}}}}
\Big(
\int_{\Omega^{^{_{l}}}}
l \gamma V^{^{_{N,l}}}_j(s,\cdot)  h^{^{_{l}}}
\,d \pi^{^{_{l}}}
-
D^{^{_{l}}}
\big(
\sqrt{h^{^{_{l}}}}
\,\big)
\Big).
\end{eqnarray}
Next we use (\ref{dirsum}) and again the variational formula (\ref{varformeq})
\begin{eqnarray}
\label{harmadik}
&&
\sup_{h^{^{_{l}}}}
\Big(
\int_{\Omega^{^{_{l}}}}
l \gamma V^{^{_{N,l}}}_j(s,\cdot)  h^{^{_{l}}}
\,d \pi^{^{_{l}}}
-
D^{^{_{l}}}
\big(
\sqrt{h^{^{_{l}}}}
\,\big)
\Big)
\\
\notag
&&
\phantom{MMM}
=
\sup_{h^{^{_{l}}}}
\sum_k w^l_k
\Big(
\int_{\Omega^{^{_{l}}}_k}
l \gamma V^{^{_{N,l}}}_j(s,\cdot)  h^{^{_{l}}}_k
\,d \pi^{^{_{l}}}_k
-
D^{^{_{l}}}_k
\big(
\sqrt{h^{^{_{l}}}_k}
\,\big)
\Big)
\\
\notag
&&
\phantom{MMM}
=
\sup_{w^{^{_{l}}}_{\cdot}}
\sum_k w^l_k
\sup_{h^{^{_{l}}}_k}
\Big(
\int_{\Omega^{^{_{l}}}_k}
l \gamma V^{^{_{N,l}}}_j(s,\cdot)  h^{^{_{l}}}_k
\,d \pi^{^{_{l}}}_k
-
D^{^{_{l}}}_k
\big(
\sqrt{h^{^{_{l}}}_k}
\,\big)
\Big)
\\
\notag
&&
\phantom{MMM}
=
\sup_{w^{^{_{l}}}_{\cdot}}
\sum_k w^l_k
\osigma
\big(
L^{^{_{l}}}_k + l \gamma V^{^{_{N,l}}}_j(s,\cdot)
\big)
\\
\notag
&&
\phantom{MMM}
=
\sup_{w^{^{_{l}}}_{\cdot}}
\sum_k w^l_k
\Big(
l \gamma \theta^{^{_{N}}}_x (s,j/N)
\expect^{^{_{l}}}_k (V^{^{_{l}}})
\\
\notag
&&
\phantom{MMMMMMMMMM}
+
\osigma
\big(
L^{^{_{l}}}_k +
l \gamma\theta^{^{_{N}}}_x (s,j/N)
\big( V^{^{_{l}}}
        -
      \expect^{^{_{l}}}_k ( V^{^{_{l}}} )
\big)
\big)
\Big)
\end{eqnarray}
In the first step we used (\ref{dirsum}). The second step is a straightforward
identity. In the third step we have used again (\ref{varformeq}) and we
introduced the notation $L^{^{_{l}}}_k$ for the infinitesimal generator of the
process restricted to $\Omega^{^{_{l}}}_k$. Finally, in the last step we use
the notation introduced at the beginning of the present paragraph.

\subsubsection{Spectral estimates}

The rest of the proof of the one block estimate relies on the
following three steps: (1) a straightforward estimate of
$\expect^{^{_{l}}}_k\big( V^{^{_{l}}} \big)$ and
$\var^{^{_{l}}}_k\big( V^{^{_{l}}} \big)$; (2) a lower bound  of
order $\sim l^{-2}$ on the spectral gap of $L^{^{_{l}}}_k$, valid
uniformly in $k\in\Z$; (3) combining these two and the
perturbational bound (\ref{gappertbound}), an upper bound on
$\osigma(\dots)$ appearing in the last expression.

\begin{lemma}
\label{expvar} There exist  constant $C(K)<\infty$ for every
$K>K_0$, such that for any $l$ and $k$ the following bounds hold:
\begin{equation}
\label{expvarbound} \left|\expect^{^{_{l}}}_k\big( V^{^{_{l}}}
\big)\right| \le C(K) l^{-1}, \qquad \var^{^{_{l}}}_k\big(
V^{^{_{l}}} \big) \le C(K) l^{-1}.
\end{equation}
\end{lemma}

\begin{proof}
For $|k|>Kl$ there is nothing to prove, so let $|k|\le Kl$.
Restricted on $\Omega^{^{_{l}}}_k$
\begin{equation*}
V^{^{_{l}}}
=
\Psi^{^{_{l}}}-\widehat\Psi(k/l)
-
\big( \Psi^{^{_{l}}}-\widehat\Psi(k/l) \big) \ind_{\{\alpha
|\Psi^{^{_{l}}}|>K\}}.
\end{equation*}
Consequently,
\begin{eqnarray*}
&&
\left|
\expect^{^{_{l}}}_k \big( V^{^{_{l}}} \big)
\right|
\le
2\left|  \expect^{^{_{l}}}_k \big(
\Psi^{^{_{l}}}-\widehat\Psi(k/l) \big) \right| +
\expect^{^{_{l}}}_k \big( \big|\Psi^{^{_{l}}}-\expect^{^{_{l}}}_k
 \Psi^{^{_{l}}}\big| \ind_{\{\alpha |\Psi^{^{_{l}}}|>K\}}\big).
\end{eqnarray*}
By the equivalence of ensembles (see e.g. Appendix 2. of
\cite{kipnislandim} and also \cite{landimsethuramanvaradhan})
\[
\left|  \expect^{^{_{l}}}_k \big( \Psi^{^{_{l}}}-\widehat\Psi(k/l)
\big) \right|\leq C(K) l^{-1}.
\]
The second term can be estimated with the Cauchy-Schwarz
inequality and with large deviation techniques (noting that
because of the growth conditions on the rates we can choose such
$\alpha>0$ that $\alpha^{-1} K>\left|\expect^{^{_{l}}}_k
\Psi^{^{_{l}}}\right|$ uniformly for $|k|<K l$) and it can be
easily shown to be smaller order then the first one.
$\var^{^{_{l}}}_k\big( V^{^{_{l}}} \big)$ may be estimated with
similar methods.
\end{proof}

\begin{lemma}
\label{gap}
There exists a constant $C<\infty$, independent of $l$ and $k$,  such that
for any $f\in L^2(\Omega^{^{_{l}}}_k, \pi^{^{_{l}}}_k)$
\begin{equation}
\label{gapbound}
\var^{^{_{l}}}_k\big( f \big)
\le
C l^{2} D^{^{_{l}}}_k\big( f \big).
\end{equation}
\end{lemma}

\begin{proof}
For the details of the proof of this gap-estimate we refer to
\cite{luyau}, \cite{landimsethuramanvaradhan},
\cite{kipnislandim}. 
For models
with bounded $z$-variable, $-\infty<\zmin<\zmax<\infty$, 
we note that 
\begin{equation*}
c(x,y)\ge \alpha\, r(x) \ind_{\{x>\zmin, y<\zmax\}}.
\end{equation*}
with some positive constant $\alpha$.
Thus, it is sufficient to prove the gap estimate for the 
\emph{reversible} process with rates 
$\tilde c(x,y):=r(x) \ind_{\{x>\zmin, y<\zmax\}}$, which has the
same ergodic stationary measures $\pi^l_k$ as our original process.
For this latter process the induction
steps of \cite{landimsethuramanvaradhan} 
or Appendix 3 of \cite{kipnislandim} 
apply without any essential modification.

In  \cite{landimsethuramanvaradhan} the statement is proved for zero range
model with rate function satisfying condition (D). Minor formal (but not
essential) modifications of that argument yield the result for the
bricklayers' models with rate functions satisfying condition (D). 
\end{proof}

\noindent
{\bf Remark:} 
Actually we could consider a wider class of models with unbounded
spin space, by imposing
\begin{equation*}
\inf_{y}c(x,y)\ge \alpha \,r(x) 
\end{equation*}
with some positive constant $\alpha$ and $r(x)$ obeying condition
(D).

We remark that there exists a constant $C$ depending only on the solution
$u(t,x)$ of the Burgers' equation (\ref{burgers}), and another constant
$C(K)$ which depends also on the cutoff level $K$, such that
\begin{eqnarray}
\label{thetabound}
\sup_{\substack{0<t<T\\ j\in \T^N}}
\big|
\theta^{^{_N}}_x (t,j/N)
\big|
&\leq&
C,
\\[5pt]
\label{supbound}
{\big\Vert
V^{l}-\expect^{^{_{l}}}_{k} V^{l}
\big\Vert}_\infty
&\leq&
C(K)
\end{eqnarray}

Now, combining (\ref{gappertbound}), (\ref{expvarbound}),  (\ref{gapbound}),
(\ref{thetabound}) and (\ref{supbound}), we get the following upper bound,
which holds for every sufficiently small $\gamma$:%
\begin{equation*}
\osigma
\big(
L^{^{_{l}}}_k
+
l \gamma\theta^{^{_{N}}}_x (s,j/N)
\big( V^{^{_{l}}}
        -
      \expect^{^{_{l}}}_k ( V^{^{_{l}}})
\big)
\big)
\le
\frac{C_1(K) l^3 \gamma^2} {1-C_2(K) \gamma l^3}
\end{equation*}
Setting
\begin{equation*}
\gamma:=\gamma_0 l^{-3} \qquad \mathrm{with} \quad \gamma_0
<\min\big\{1,\big(2  C_2(K)\big)^{-1}\big\}
\end{equation*}
we have
\begin{equation*}
\osigma
\big(
L^{^{_{l}}}_k
+
l \gamma\theta^{^{_{N}}}_x (s,j/N)
\big( V^{^{_{l}}}
        -
      \expect^{^{_{l}}}_k ( V^{^{_{l}}} )
\big)
\big)
\le
C(K)\gamma_0^2 l^{-3}.
\end{equation*}
Collecting all the estimates and going backwards through
(\ref{harmadik}), (\ref{masodik}), (\ref{elso}), we find eventually
\begin{equation*}
\log \expect_{\pi^{^{_{N}}}}
\exp
\Big\{
\int_0^{N^{1+\beta}t}
\gamma V^{^{_{N,l}}}\big(s,\uzeta^{^{_{N}}} (s)\big)\,ds
\Big\}
\le
C(K) \gamma_0 N^{2+\beta}l^{-4}
\end{equation*}
and
\begin{equation}
\label{abound} \int_0^t A^{^{_{N,l}}}_K (s) \,ds \leq
C(K)(N^{-3\beta}l^3 + N l^{-1})
\end{equation}
Consequently, from (\ref{abound}), (\ref{bbound}) and (\ref{m=a+b}),  with any
fixed $K<\infty$ we have
\begin{equation}
\label{oneblockerror2} \int_0^tM^{^{_{N,l}}}(s)\,ds \leq
\epsilon(K) {\cal O} \big( N^{1-2 \beta} \big) +
C(K)(N^{-3\beta}l^3 + N l^{-1})
\end{equation}
where $C(K)$ is a finite constant which may increase to infinity
as $K\to\infty$, and $\epsilon(K)\to 0$ as $K\to\infty$.

\subsection{End of proof}
\label{konyec}

We put together (\ref{sumup3}) and (\ref{oneblockerror2}) to get,
for any $K<\infty$ fixed (with a $C$ not depending on $K$)
\begin{eqnarray*}
&&
H^{^{_{N}}}(t)
\le
\phantom{+} H^{^{_{N}}}(0) + C \int_0^t H^{^{_{N}}}(s)ds +
\epsilon(K) {\cal O} \big( N^{1-2 \beta} \big)
\\
\notag
&& \qquad\qquad
+
{\cal O}
\big(
N^{1-3\beta}
\vee
N^{-\beta}l
\vee
N l^{-1}
\vee N^{-3\beta}l^{3}
\big).
\end{eqnarray*}
If
\begin{equation*}
0<\beta<\frac15
\end{equation*}
then we can choose
\begin{equation*}
N^{2\beta}\ll l \ll N^{(1+\beta)/3}
\end{equation*}
which ensures
\begin{equation*}
{\cal O}
\big(
N^{1-3\beta}
\vee
N^{-\beta} l
\vee
N l^{-1}
\vee
N^{-3\beta}l^{3}
\big)
=
o\big( N^{1-2\beta} \big).
\end{equation*}
Thus for every $K<\infty$%
\begin{eqnarray*}
H^{^{_{N}}}(t)
\le
H^{^{_{N}}}(0) + C \int_0^t H^{^{_{N}}}(s)ds + \epsilon(K)
N^{1-2 \beta}+o\big( N^{1-2\beta} \big),
\end{eqnarray*}
where $\lim_{K\rightarrow \infty}\epsilon(K)=0$, and from
Gromwall indeed (\ref{main}) follows, uniformly for $t\in[0,T]$.


\vskip1cm
\noindent
{\Large\bf Acknowledgment:}
It is a pleasure of the authors to thank J\'ozsef Fritz many illuminating
discussions on the topics of hydrodynamical limits in general and on some
particular technical aspects of the present work. B.T. also thanks the kind
hospitality of Institut Henri Poincar\'e, where part of this work was done.


\vskip1cm

\hbox{\sc
\vbox{\noindent
\hsize66mm
B\'alint T\'oth\\
Institute of Mathematics\\
Technical University Budapest\\
Egry J\'oszef u. 1.\\
H-1111 Budapest, Hungary\\
{\tt balint@math.bme.hu}
}
\hskip5mm
\vbox{\noindent
\hsize66mm
Benedek Valk\'o\\
Institute of Mathematics\\
Technical University Budapest\\
Egry J\'oszef u. 1.\\
H-1111 Budapest, Hungary\\
{\tt valko@math.bme.hu}
}
}


\begin{thebibliography}{99}

\bibitem{balazs1}
M. Bal\'azs:
Microscopic structure of the shock in a domain growth model.
{\sl Journal of Statistical Physics} {\bf 105} no. 1-2 (2001)

\bibitem{balazs2}
M. Bal\'azs:
Growth fluctuations in interface models.
{\sl preprint} (2001)

\bibitem{cocozza}
C. Cocozza:
Processus des misanthropes.
{\sl Zeitschrift f\"ur Wahrscheinlichkeitstheorie und verwandte Gebiete}
{\bf 70}: 509-523 (1985)

\bibitem{espositomarrayau}
R. Esposito, R. Marra, H.T. Yau:
Diffusive limit of asymmetric simple exclusion.
{\sl Reviews of Mathematical Physics} {\bf 6}: 1233-1267 (1994)


\bibitem{fritz}
J. Fritz:
{\sl An Introduction to the Theory of Hydrodynamic Limits.\/}
Lectures in Mathematical Sciences.
Graduate School of Mathematics, Univ. Tokyo, 2001.


\bibitem{kipnislandim}
C. Kipnis, C. Landim:
{\sl Scaling Limits of Interacting Particle Systems.\/}
Springer, 1999.

\bibitem{komoriya}
K. Komoriya:
Hydrodynamic limit for asymmetric mean zero exclusion processes with speed
change.
{\sl Annales de l'Institut Henri Poincar\'e --- Probabilit\'es et
Statistiques} {\bf 34}: 767-797 (1998)


\bibitem{landimsethuramanvaradhan}
C. Landim, S. Sethuraman, S.R.S. Varadhan:
Spectral gap for zero range dynamics.
{\sl Annals of Probability} {\bf 24}: 1871-1902 (1986)

\bibitem{luyau}
S.L. Lu, H-T. Yau:
Spectral gap and logarothmic Sobolev inequality for Kawasaki and
Glauber dynamics.
{\sl Communications in Mathematical Physics} {\bf 156}: 399-433 (1993)

\bibitem{rezakhanlou2}
F. Rezakhanlou:
Microscopic structure of shocks in one conservation laws.
{\sl Annales de l'Institut Henri Poincar\'e --- Analyse Non Lineaire}
{\bf 12}: 119-153 (1995)

\bibitem{seppalainen}
T. Sepp\"al\"ainen:
Perturbation of the equilibrium for a  totally asymmetric stick
process in one dimension.
{\sl Annals of Probability} {\bf 29}: 176-204 (2001)

\bibitem{yau}
H.T. Yau:
Relative entropy and hydrodynamics of Ginzburg-Landau models.
{\sl Lett. Math. Phys.\/} {\bf 22}: 63-80 (1991)



\end{thebibliography}
\end{document}